\documentclass[12pt]{article}
\usepackage{amsfonts}
\usepackage{full page,
amssymb, amscd, amsmath,graphicx, 
}

 \title{{\bf On functors between module categories for associative algebras
and for $\N$-graded vertex algebras}}
 \author{Yi-Zhi Huang and Jinwei Yang}
    \date{}
    \begin{document}
    \bibliographystyle{alpha}

\newtheorem{thm}{Theorem}[section]
\newtheorem{defn}[thm]{Definition}
\newtheorem{prop}[thm]{Proposition}
\newtheorem{cor}[thm]{Corollary}
\newtheorem{lemma}[thm]{Lemma}
\newtheorem{rema}[thm]{Remark}
\newtheorem{app}[thm]{Application}
\newtheorem{prob}[thm]{Problem}
\newtheorem{conv}[thm]{Convention}
\newtheorem{conj}[thm]{Conjecture}
\newtheorem{cond}[thm]{Condition}
    \newtheorem{exam}[thm]{Example}
\newtheorem{assum}[thm]{Assumption}
     \newtheorem{nota}[thm]{Notation}
\newcommand{\halmos}{\rule{1ex}{1.4ex}}
\newcommand{\pfbox}{\hspace*{\fill}\mbox{$\halmos$}}

\newcommand{\nn}{\nonumber \\}

 \newcommand{\res}{\mbox{\rm Res}}
 \newcommand{\ord}{\mbox{\rm ord}}
\renewcommand{\hom}{\mbox{\rm Hom}}
\newcommand{\edo}{\mbox{\rm End}\ }
 \newcommand{\pf}{{\it Proof.}\hspace{2ex}}
 \newcommand{\epf}{\hspace*{\fill}\mbox{$\halmos$}}
 \newcommand{\epfv}{\hspace*{\fill}\mbox{$\halmos$}\vspace{1em}}
 \newcommand{\epfe}{\hspace{2em}\halmos}
\newcommand{\nord}{\mbox{\scriptsize ${\circ\atop\circ}$}}
\newcommand{\wt}{\mbox{\rm wt}\ }
\newcommand{\swt}{\mbox{\rm {\scriptsize wt}}\ }
\newcommand{\lwt}{\mbox{\rm wt}^{L}\;}
\newcommand{\rwt}{\mbox{\rm wt}^{R}\;}
\newcommand{\slwt}{\mbox{\rm {\scriptsize wt}}^{L}\,}
\newcommand{\srwt}{\mbox{\rm {\scriptsize wt}}^{R}\,}
\newcommand{\clr}{\mbox{\rm clr}\ }
\newcommand{\tr}{\mbox{\rm Tr}}
\newcommand{\C}{\mathbb{C}}
\newcommand{\Z}{\mathbb{Z}}
\newcommand{\R}{\mathbb{R}}
\newcommand{\Q}{\mathbb{Q}}
\newcommand{\N}{\mathbb{N}}
\newcommand{\CN}{\mathcal{N}}
\newcommand{\F}{\mathcal{F}}
\newcommand{\I}{\mathcal{I}}
\newcommand{\V}{\mathcal{V}}
\newcommand{\one}{\mathbf{1}}
\newcommand{\BY}{\mathbb{Y}}
\newcommand{\ds}{\displaystyle}

        \newcommand{\ba}{\begin{array}}
        \newcommand{\ea}{\end{array}}
        \newcommand{\be}{\begin{equation}}
        \newcommand{\ee}{\end{equation}}
        \newcommand{\bea}{\begin{eqnarray}}
        \newcommand{\eea}{\end{eqnarray}}
         \newcommand{\lbar}{\bigg\vert}
        \newcommand{\p}{\partial}
        \newcommand{\dps}{\displaystyle}
        \newcommand{\bra}{\langle}
        \newcommand{\ket}{\rangle}

        \newcommand{\ob}{{\rm ob}\,}
        \renewcommand{\hom}{{\rm Hom}}

\newcommand{\A}{\mathcal{A}}
\newcommand{\Y}{\mathcal{Y}}
    \maketitle

\newcommand{\dlt}[3]{#1 ^{-1}\delta \bigg( \frac{#2 #3 }{#1 }\bigg) }

\newcommand{\dlti}[3]{#1 \delta \bigg( \frac{#2 #3 }{#1 ^{-1}}\bigg) }

 \makeatletter
\newlength{\@pxlwd} \newlength{\@rulewd} \newlength{\@pxlht}
\catcode`.=\active \catcode`B=\active \catcode`:=\active
\catcode`|=\active
\def\sprite#1(#2,#3)[#4,#5]{
   \edef\@sprbox{\expandafter\@cdr\string#1\@nil @box}
   \expandafter\newsavebox\csname\@sprbox\endcsname
   \edef#1{\expandafter\usebox\csname\@sprbox\endcsname}
   \expandafter\setbox\csname\@sprbox\endcsname =\hbox\bgroup
   \vbox\bgroup
  \catcode`.=\active\catcode`B=\active\catcode`:=\active\catcode`|=\active
      \@pxlwd=#4 \divide\@pxlwd by #3 \@rulewd=\@pxlwd
      \@pxlht=#5 \divide\@pxlht by #2
      \def .{\hskip \@pxlwd \ignorespaces}
      \def B{\@ifnextchar B{\advance\@rulewd by \@pxlwd}{\vrule
         height \@pxlht width \@rulewd depth 0 pt \@rulewd=\@pxlwd}}
      \def :{\hbox\bgroup\vrule height \@pxlht width 0pt depth
0pt\ignorespaces}
      \def |{\vrule height \@pxlht width 0pt depth 0pt\egroup
         \prevdepth= -1000 pt}
   }
\def\endsprite{\egroup\egroup}
\catcode`.=12 \catcode`B=11 \catcode`:=12 \catcode`|=12\relax
\makeatother

\def\hboxtr{\FormOfHboxtr} 
\sprite{\FormOfHboxtr}(25,25)[0.5 em, 1.2 ex] 

:BBBBBBBBBBBBBBBBBBBBBBBBB | :BB......................B |
:B.B.....................B | :B..B....................B |
:B...B...................B | :B....B..................B |
:B.....B.................B | :B......B................B |
:B.......B...............B | :B........B..............B |
:B.........B.............B | :B..........B............B |
:B...........B...........B | :B............B..........B |
:B.............B.........B | :B..............B........B |
:B...............B.......B | :B................B......B |
:B.................B.....B | :B..................B....B |
:B...................B...B | :B....................B..B |
:B.....................B.B | :B......................BB |
:BBBBBBBBBBBBBBBBBBBBBBBBB |

\endsprite
\def\shboxtr{\FormOfShboxtr} 
\sprite{\FormOfShboxtr}(25,25)[0.3 em, 0.72 ex] 

:BBBBBBBBBBBBBBBBBBBBBBBBB | :BB......................B |
:B.B.....................B | :B..B....................B |
:B...B...................B | :B....B..................B |
:B.....B.................B | :B......B................B |
:B.......B...............B | :B........B..............B |
:B.........B.............B | :B..........B............B |
:B...........B...........B | :B............B..........B |
:B.............B.........B | :B..............B........B |
:B...............B.......B | :B................B......B |
:B.................B.....B | :B..................B....B |
:B...................B...B | :B....................B..B |
:B.....................B.B | :B......................BB |
:BBBBBBBBBBBBBBBBBBBBBBBBB |

\endsprite

\vspace{2em}

\begin{abstract}
We prove that the weak associativity for modules for vertex algebras
are equivalent to a residue formula for iterates of vertex operators,
obtained using the weak associativity
and the lower truncation property
of vertex operators, together with a known
formula expressing products of components of vertex operators
as linear combinations of iterates of components of vertex operators.
By requiring that these two formulas instead of the
commutator formula hold, we construct a functor $S$
from the category of modules for Zhu's algebra of a vertex operator algebra $V$
to the category of $\N$-gradable weak $V$-modules.
We prove that $S$ has a universal property and
the functor $T$ of taking top levels of $\N$-gradable weak $V$-modules
is a left inverse of $S$. In particular, $S$
is equal to a functor implicitly given by Zhu and explicitly 
constructed by Dong, Li and Mason
and we obtain a new construction without using relations corresponding to
the commutator formula.
The hard part of this new construction
is a technical theorem stating roughly that
in a module for Zhu's algebra, the relation corresponding to the residue formula
mentioned above can in fact be obtained from the relations corresponding to
the action of Zhu's algebra.
\end{abstract}



\renewcommand{\theequation}{\thesection.\arabic{equation}}
\renewcommand{\thethm}{\thesection.\arabic{thm}}
\setcounter{equation}{0} \setcounter{thm}{0} 
\date{}
\maketitle

\setcounter{equation}{0}

\section{Introduction}

In the study of vertex operator algebras and their modules,
commutator formulas for vertex operators often play an important role.
Because of the commutator formula, the components of
the vertex operators of a vertex operator algebra
form a Lie algebra. The vertex operator algebra itself
and its modules become modules for this Lie algebra.
This fact allows one to apply many techniques in the representation
theory of Lie algebras to  the study of vertex
operator algebras and their modules. For example,
the proof by Zhu \cite{Z} of
the modular invariance of the space spanned by the $q$-traces
of products of vertex operators
depends heavily on the commutator formula, especially in the
case of $q$-traces of products of at least two vertex operators.

On the other hand, the main interesting mathematical objects in
the representation theory of vertex operator algebras are intertwining
operators, not just vertex operators for vertex operator algebras
and for modules. For intertwining operators, in general there is no
commutator formula and therefore the Lie-algebra-theoretic techniques
cannot be applied. Thus it is necessary
to develop different, non-Lie-algebra-theoretic method to study intertwining
operators. One example is the modular invariance of the space spanned by the $q$-traces
of products of at least two intertwining operators. This
modular invariance was conjectured first by physicists and was
used explicitly by Moore and Seiberg in \cite{MS} to derive their
polynomial equations. But its proof was given by the first author in the
preprint of \cite{H2} 13 years after
Zhu's proof of his modular invariance result in \cite{Z}. The difficulty
lies mainly in the development of new method without using the commutator
formula on which Zhu's proof depends heavily.
Instead of any commutator formula (which in fact does not exist
in general), the first author
used the associativity for intertwining operators proved in \cite{H1}.

Because of the limitation of the commutator formula
discussed above, it is desirable to prove results in the representation theory
of vertex operator algebras, even only for modules, without using
the commutator formula so that it will be relatively easy to see whether such results
can be generalized to intertwining operators.

Note that for modules for vertex operator algebras, the weak associativity
can be taken to be the main axiom. It is natural to expect that
all the results for modules can be proved using the weak associativity
without the commutator formula. However, in many problems in the representation
theory of vertex operator algebras and its applications,
one has to work with components of vertex operators instead of vertex operators
directly. Thus we need formulas expressed using residues of
formal series obtained from products and iterates of vertex operators
that are equivalent to the weak associativity.

It is known from \cite{DLM} and \cite{L} (see also \cite{LL}) that for a vertex algebra,
the weak associativity for vertex operators gives a
formula that expresses products of components of vertex operators
as linear combinations of iterates of components of vertex operators.
Unfortunately, this formula is not equivalent to the weak associativity.
In the present paper, we give a formula for the residues of certain formal
series involving iterates of vertex operators
obtained using the weak associativity
and the lower truncation property of vertex operators. Then we
prove that the weak associativity is
equivalent to this formula together with the formula given in \cite{DLM} and \cite{L}.
Thus, in principle, all the results for modules for vertex operator algebras
can be proved using these two formulas without the help of the commutator formula.

We apply this result to give a new construction of $\N$-gradable weak modules
for an $\N$-graded vertex
algebra $V$ from modules for Zhu's algebra $A(V)$. We construct a functor $S$ from
the category of  $A(V)$-modules to the category of $\N$-gradable weak
$V$-modules. We prove that $S$ has a universal property and
the functor $T$ of taking top levels of $\N$-gradable weak $V$-modules
is a left inverse of $S$. In particular, $S$
is equal to the functor constructed by Dong, Li and Mason in \cite{DLM}
in the case $n=0$
and we achieve our goal of finding a new construction
without using relations corresponding to
the commutator formula. The hard part of this new construction
is a technical theorem stating roughly that
in a module for Zhu's algebra, the relation corresponding to the residue formula
given in the present paper (see above) can in fact be obtained from the relations
corresponding to the action of Zhu's algebra.

The importance of our construction is that it allows us to use the method in
the present paper  to
give constructions and prove results in the cases that there is no commutator formula.
The main example that motivated the present paper is a construction
of generalized twisted modules associated to an infinite order automorphism of
a vertex operator algebra in the sense of the first author in \cite{H3}.
For these generalized twisted modules, the twisted
vertex operators involve the logarithm of the variable and thus do not have a
commutator formula. We shall discuss this construction in a future publication.

This paper is organized as follows: We recall and derive some
residue formulas for vertex operators from the weak associativity
in Section 2. We prove that these residue formulas are equivalent to
the weak associativity in Section 3. In Section 4, we construct our functor $S$.
In Section 5, we prove a university property of $S$ and prove that the
functor $T$ mentioned above is a left inverse of $S$.

\paragraph{Acknowledgments}
The results in this paper was presented by the second author 
in the conference ``Representation theory XIII,''  Dubrovnik, Croatia, 
June 21-27, 2013.
The second author would like to thank Drazen Adamovic and 
Antun Milas for inviting him to give the talk in the conference.
This research is supported in part by NSF
grant PHY-0901237.

\section{Residue formulas from the weak associativity}

In this section, we recall and derive some consequences of
the weak associativity.

Let $(V, Y, \one)$ be a vertex algebra and $(W, Y_{W})$ a $V$-module.
(Note that since $V$ as a vertex algebra might not have a
grading, a $V$-module also might not have a grading.)
The weak associativity for $W$ states that for $u, v\in V$, $w\in W$ and
$l$ a nonnegative integer such that $u_{n}w=0$ for $n\ge l$,
we have
$$(x_{0}+x_{2})^{l}Y_{W}(u, x_{0}+x_{2})Y_{W}(v, x_{2})w
=(x_{0}+x_{2})^{l}Y_{W}(Y(u, x_{0})v, x_{2})w.$$
We recall the following result from
\cite{DLM} and \cite{L} (Proposition 4.5.7 in \cite{LL}):

\begin{prop}\label{p2}
Let $W$ be a $V$-module
and let $u, v \in V$, $p, q \in \Z$ and $w \in W$. 
Let $l$ be a nonnegative integer such that
\[
u_nw = 0\;\;\; {\rm for}\; n \geq l
\]
and let $k$ be a nonnegative integer such that
$k-q$ is positive and
\[
v_nw = 0\;\;\;{\rm for}\; n \geq k.
\]
Then
\begin{equation}\label{eq1}
u_pv_qw = \sum_{i = 0}^{k-q-1}\sum_{j = 0}^l\binom{p - l}{i}\binom{l}{j}
(u_{p-l-i+j}v)_{q+l+i-j}w.
\end{equation}\epfv
\end{prop}

\begin{rema}\label{eq1-rmk}
{\rm Let
$$p(x_{0}, x_{2})=\sum_{i=0}^{k-q-1}\binom{p-l}{i}x_{0}^{p-l-i}x_{2}^{i}.$$
Then the formula (\ref{eq1}), the formula
\begin{eqnarray}\label{eq1-0}
\lefteqn{\res_{x_{1}}\res_{x_{2}}x_{1}^{p}x_{2}^{q}Y(u, x_{1})
Y(v, x_{2})w}\nn
&&=\res_{x_{0}}\res_{x_{2}}p(x_{0}, x_{2})x_{2}^{q}(x_{0}+x_{2})^{l}
Y(Y(u, x_{0})v, x_{2})w
\end{eqnarray}
and the formula
\begin{eqnarray}\label{eq1-1}
\lefteqn{\res_{x_{0}}\res_{x_{2}}(x_{0}+x_{2})^{p}x_{2}^{q}Y(u, x_{0}+x_{2})
Y(v, x_{2})w}\nn
&&=\res_{x_{0}}\res_{x_{2}}p(x_{0}, x_{2})x_{2}^{q}(x_{0}+x_{2})^{l}
Y(Y(u, x_{0})v, x_{2})w
\end{eqnarray}
are equivalent to each other.
See the proof of Proposition 4.5.7 in \cite{LL}.}
\end{rema}

From the lower truncation property for vertex operators, we have the following:
\begin{prop}
For $W$, $u, v$, $w$, $l, k$ as in Proposition \ref{p2}, we have
\begin{eqnarray}
\res_{x_{2}}x_{2}^{q}(x_{0}+x_{2})^{l}Y_{W}(u, x_{0}+x_{2})
Y_{W}(v, x_{2})w&=&0,\label{upper-truncation-1}\\
\res_{x_{2}}x_{2}^{q}(x_{0}+x_{2})^{l}Y_{W}(Y(u, x_{0})v, x_{2})
w&=&0\label{upper-truncation-2}
\end{eqnarray}
for $q\ge k$.
\end{prop}
\pf
Note that $x_{2}^{q}(x_{0}+x_{2})^{l}Y_{W}(u, x_{0}+x_{2})$ contains only terms
with the powers of $x_{2}$ larger than or equal to $q\ge k$. Hence
the left hand side of (\ref{upper-truncation-1}) are linear combinations of
$u_{m}v_{n}w$ for $m\in \Z$ and $n\ge k$ with coefficients in $\C((x_{0}))$.
Since $v_{n}w=0$ for $n\ge k$, (\ref{upper-truncation-1}) holds.

By the weak associativity, the left-hand side of
(\ref{upper-truncation-2}) is equal to the left-hand side of (\ref{upper-truncation-1}).
Thus (\ref{upper-truncation-2}) also holds.
\epf

\begin{rema}\label{upper-truncation-rmk}
{\rm Note that the proof of (\ref{upper-truncation-1}) uses only the property
of $Y_{W}$ that $Y_{W}(v, x)w\in W((x))$ for $v\in V$ and $w\in W$. Thus for
a vector space $W$ and a linear map $Y_W$ from $V\otimes W$ to $W((x))$,
(\ref{upper-truncation-1}) still holds. On the other hand,
(\ref{upper-truncation-2}) does not hold in general for such a vector space
$W$ and such a map $Y_W$.}
\end{rema}

From (\ref{upper-truncation-2}), we obtain immediately:

\begin{cor}
For $W$, $u, v$, $p, q$, $w$, $l, k$ as in Proposition \ref{p2} except that
$k-q$ does not have to be positive, we have
\begin{equation}\label{upper-truncation-3}
\res_{x_{0}}\res_{x_{2}}x_{0}^{p-l-i}x_{2}^{q+i}(x_{0}+x_{2})^{l}Y_{W}(Y(u, x_{0})v, x_{2})
w=0
\end{equation}
for $i\ge k-q$.\epf
\end{cor}

\begin{rema}\label{upper-truncation-rmk-2}
{\rm Assume that $V$ is a $\Z$-graded vertex algebra
and $W$ is an $\N$-gradable weak $V$-module.
Here by  a  $\Z$-graded vertex algebra we mean a
vertex algebra
$V=\coprod_{n\in \Z}V_{(n)}$ such that
for $u\in V$,
$$[d, Y(u, x)]=x\frac{d}{dx}Y(u, x)+Y(du, x),$$
where $d: V\to V$ is defined by $du=nu$ for $u\in V_{(n)}$.
For $u\in V_{(n)}$, we use $\wt u$ to denote $n$.
By an $\N$-gradable weak $V$-module, we mean a module $W=\coprod_{n\in \N}W_{(n)}$
for $V$ when $V$ is  viewed as a vertex algebra such that for $u\in V_{(m)}$
$u_{n}=\res_{x}x^{n}Y_{W}(u, x)$ maps $W_{(k)}$ to $W_{(k+m-n-1)}$.
For $w\in W_{(n)}$, we use $\deg w$ to denote $n$.
When $u$, $v$ and $w$ are homogeneous,
we take $k=\wt v+\deg w$ and $l=\wt u+\deg w$ and
set $K = p + q + 2-\wt u-\wt v$ and $m = p - l - i$. Then
(\ref{upper-truncation-3}) becomes
\begin{equation}\label{upper-truncation-4}
\res_{x_{0}}\res_{x_{2}}x_{0}^{m}x_{2}^{K+\swt v-m-\deg w-2}
(x_{0}+x_{2})^{\swt u+\deg w}
Y_{W}(Y(u, x_{0})v, x_{2})
w=0
\end{equation}
for $m \leq M-2\deg w - 2$. The formula (\ref{upper-truncation-4})
can be further written without using the weights of $u$ and $v$ as
\begin{equation}\label{upper-truncation-4.1}
\res_{x_{0}}\res_{x_{2}}x_{0}^{m}x_{2}^{K-m-\deg w-2}(x_{0}+x_{2})^{\deg w}
Y_{W}(Y((x_{0}+x_{2})^{L(0)}u, x_{0})x_{2}^{L(0)}v, x_{2})
w=0,
\end{equation}
which is in fact holds for general $u$ and $v$, not necessarily homogeneous.
The component form of
(\ref{upper-truncation-4}) is
\begin{equation}
\sum_{j=0}^{\swt u+\deg w}\binom{\swt u+\deg w}{ j}
(u_{j+m}v)_{K+\swt u+\swt v - j - m -2}w=0
\end{equation}
for $m \leq K-2\deg w - 2$. }
\end{rema}

\setcounter{equation}{0}
\section{The equivalence between the residue formulas and the weak associativity}

In this section, we prove that (\ref{eq1-1}) and (\ref{upper-truncation-3})
imply, and therefore are equivalent to, the weak associativity.

\begin{thm}\label{main proposition}
Let $V$ be a vertex operator algebra, $W$ a vector space and
$Y_W$ a linear map from $V\otimes W$ to $W((x))$. For $v\in V$ and $w\in W$,
as for a $V$-module, we denote the image of $u\otimes w$ under $Y_W$
by $Y_W(u, x)w$ and $\res_{x}x^{n}Y_W(u, x)w$ by $u_{n}$.
Let $u, v\in V$ and $w\in W$ and let  $k, l \in \Z$ such that
\begin{equation}\label{eq6}
v_nw = 0\;\;\; {\rm for}\; n \geq k
\end{equation}
and
\begin{equation}\label{eq7}
u_nw = 0\;\;\;{\rm for}\; n \geq l.
\end{equation}
Assume that (i) the formulas (\ref{eq1-1}) (or equivalently,
(\ref{eq1-0}) or (\ref{eq1})) holds for $p, q\in \Z$ such that $q<k$
and (ii) (\ref{upper-truncation-3}) holds for for $p, q\in \Z$ such that $i\ge k-q$.
Then
\begin{equation}\label{wk-assoc}
(x_{0}+x_{2})^{l}Y_W(u, x_0 + x_2)Y_W(v, x_2)w
=(x_{0}+x_{2})^{l}Y_W(Y(u, x_0)v, x_2)w.
\end{equation}
\end{thm}
\pf
The Laurent polynomial $p(x_{0}, x_{2})$ in Remark \ref{eq1-rmk}
is in fact the first $k-q-1$ terms of the formal series $(x_{0}+x_{2})^{p-l}$.
But from (\ref{upper-truncation-3}), we obtain
\begin{equation}\label{upper-truncation-7}
\res_{x_{0}}\res_{x_{2}}((x_{0}+x_{2})^{p-l}-p(x_{0}, x_{2}))x_{2}^{q}
((x_{0}+x_{2})^{l}Y_{W}(Y(u, x_{0})v, x_{2}))
w=0.
\end{equation}
Combining (\ref{eq1-1}) and (\ref{upper-truncation-7}),
we obtain
\begin{eqnarray*}
\lefteqn{\res_{x_{0}}\res_{x_{2}}(x_{0}+x_{2})^{p-l}x_{2}^{q}((x_{0}+x_{2})^{l}
Y_{W}(u, x_{0}+x_{2})
Y_{W}(v, x_{2}))w}\nn
&&=\res_{x_{0}}\res_{x_{2}}(x_{0}+x_{2})^{p-l}x_{2}^{q}((x_{0}+x_{2})^{l}
Y_{W}(Y(u, x_{0})v, x_{2}))w,
\end{eqnarray*}
or equivalently,
\begin{equation}\label{eq1-2}
\res_{x_{0}}\res_{x_{2}}(x_{0}+x_{2})^{p-l}x_{2}^{q}((x_{0}+x_{2})^{l}
(Y(u, x_{0}+x_{2})
Y(v, x_{2})-Y_{W}(Y(u, x_{0})v, x_{2})))w=0,
\end{equation}

On the other hand, by Remark \ref{upper-truncation-rmk}, (\ref{upper-truncation-1})
holds for $m\ge k$. In particular, we have
\begin{equation}\label{upper-truncation-5}
\res_{x_{0}}\res_{x_{2}}x_{0}^{p-l-i}x_{2}^{q+i}(x_{0}+x_{2})^{l}Y_{W}(u, x_{0}+x_{2})
Y_{W}(v, x_{2})
w=0
\end{equation}
for $i\ge k-q$. From (\ref{upper-truncation-3}) and (\ref{upper-truncation-5}),
we obtain
\begin{equation}\label{upper-truncation-6}
\res_{x_{0}}\res_{x_{2}}x_{0}^{p-l-i}x_{2}^{q+i}(x_{0}+x_{2})^{l}(Y_{W}(u, x_{0}+x_{2})
Y_{W}(v, x_{2})-Y_{W}(Y(u, x_{0})v, x_{2}))
w=0
\end{equation}
for $i\ge k-q$.
Combining (\ref{eq1-2}) and (\ref{upper-truncation-6}), we obtain
\begin{eqnarray}\label{eq1-3}
\lefteqn{\res_{x_{0}}\res_{x_{2}}\sum_{i=0}^{k-q-1}\binom{p-l}{i}
x_{0}^{p-l-i}x_{2}^{q+i}\cdot}\nn
&&\qquad\qquad\qquad \cdot ((x_{0}+x_{2})^{l}(Y(u, x_{0}+x_{2})
Y(v, x_{2})-Y_{W}(Y(u, x_{0})v, x_{2})))w\nn
&&=0.
\end{eqnarray}

We now use induction on $k-q-1$ to prove
\begin{equation}\label{eq1-3-0}
\res_{x_{0}}\res_{x_{2}}
x_{0}^{p-l}x_{2}^{q} ((x_{0}+x_{2})^{l}(Y(u, x_{0}+x_{2})
Y(v, x_{2})-Y_{W}(Y(u, x_{0})v, x_{2})))w=0
\end{equation}
for $p\in \Z$ and $q<k$.
When $k-q-1=0$, (\ref{eq1-3}) becomes
$$\res_{x_{0}}\res_{x_{2}}
x_{0}^{p-l}x_{2}^{q} ((x_{0}+x_{2})^{l}(Y(u, x_{0}+x_{2})
Y(v, x_{2})-Y_{W}(Y(u, x_{0})v, x_{2})))w=0.$$
Assuming that (\ref{eq1-3-0}) holds when $0\le k-q-1<n$. When
$k-q-1=n$, $0\le k-q-i-1<n$ for $i=1, \dots, n=k-q-1$. Since $p$ is arbitrary,
we can replace $p$ by $p-i$ for any $i\in \Z$ in (\ref{eq1-3-0}).
Thus by the induction assumption,
\begin{equation}\label{eq1-3-1}
\res_{x_{0}}\res_{x_{2}}
x_{0}^{p-l-i}x_{2}^{q+i} ((x_{0}+x_{2})^{l}(Y(u, x_{0}+x_{2})
Y(v, x_{2})-Y_{W}(Y(u, x_{0})v, x_{2})))w=0
\end{equation}
for $i=1, \dots, n=k-q-1$.
From (\ref{eq1-3-1}) for $i=1, \dots, n=k-q-1$ and (\ref{eq1-3}),
we obtain
\begin{eqnarray*}
\lefteqn{\res_{x_{0}}\res_{x_{2}}
x_{0}^{p-l}x_{2}^{q}((x_{0}+x_{2})^{l}(Y(u, x_{0}+x_{2})
Y(v, x_{2})-Y_{W}(Y(u, x_{0})v, x_{2})))w}\nn
&&=\res_{x_{0}}\res_{x_{2}}\sum_{i=0}^{k-q-1}\binom{p-l}{ i}
x_{0}^{p-l-i}x_{2}^{q+i}\cdot\nn
&&\qquad\qquad\qquad \cdot ((x_{0}+x_{2})^{l}(Y(u, x_{0}+x_{2})
Y(v, x_{2})-Y_{W}(Y(u, x_{0})v, x_{2})))w\nn
&&=0,
\end{eqnarray*}
proving (\ref{eq1-3-0}) in this case.

Taking $i=0$ in  (\ref{upper-truncation-3}), we see that
(\ref{eq1-3-0}) also holds for $p\in \Z$ and $q\ge k$. Thus
(\ref{upper-truncation-3}) holds for $p, q\in \Z$. But this means that
$$((x_{0}+x_{2})^{l}(Y(u, x_{0}+x_{2})
Y(v, x_{2})-Y_{W}(Y(u, x_{0})v, x_{2})))w=0,$$
that is, (\ref{wk-assoc}) holds.
 \epfv

The following result follows immediately from Theorem \ref{main proposition}
above and Theorems 3.6.3 and 3.11.8 in \cite{LL}.

\begin{thm}\label{main theorem}
Let $V$ be a $\Z$-graded vertex algebra, $W=\coprod_{n\in \N}W_{(n)}$
an $\N$-graded vector space and
$Y_W$ a linear map from $V\otimes W$ to $W((x))$. For $v\in V$ and $w\in W$,
we denote the image of $u\otimes w$ under $Y_W$
by $Y_W(u, x)w$ and $\res_{x}x^{n}Y_W(u, x)w$ by $u_{n}$.
Assume that $u_{n}$ maps $W_{(k)}$ to
$W_{(k+m-n-1)}$ for for $u\in V_{(m)}$ and $n\in \Z$ and $Y_{W}(\one, x)=1_{W}$.
Also assume that for $u, v\in V$, $w\in W$, there exist
$k, l \in \Z$ such that
(\ref{eq6}) and (\ref{eq7}) hold, and for $p, q\in \Z$, $u, v\in V$, $w\in W$,
the formulas (\ref{eq1-0}) (or equivalently,
(\ref{eq1}) or (\ref{eq1-1})) and (\ref{upper-truncation-3}) (or equivalently,
(\ref{upper-truncation-4})) hold.
Then $(W, Y_W)$ is an $\N$-gradable weak $V$-module.\epf
\end{thm}

\setcounter{equation}{0}
\section{A functor $S$}

In this and the next section, $(V, Y, \one)$ is a $\Z$-graded vertex algebra 
such that
$V_{(n)}=0$ for $n<0$. We shall call such a $\Z$-graded vertex algebra
an {$\N$-graded vertex algebra}. In this
section, by dividing relations such that (\ref{eq1-1}) and
(\ref{upper-truncation-3}) hold,
we construct a functor $S$ from the category of modules for
Zhu's algebra $A(V)$ to the category of $\N$-gradable weak $V$-modules.

We first recall the definition of Zhu's algebra in \cite{Z}.
Let $O(V)$ be the subspace of $V$ spanned by elements of the form
\[
\res_x x^{n}Y((x+1)^{L(0)}u, x)v
\]
for $u, v \in V$ and $n\le -2$. Zhu's algebra associated to $V$
is  the quotient space $A(V)=V/O(V)$ equipped with the
multiplication
$$u * v = \res_x x^{-1}Y((x+1)^{L(0)}u, x)v$$
for $u, v \in V$.
It was proved in \cite{Z} that $A(V)$ equipped with the product $*$ is
an associative algebra.

Let  $(W, Y_{W})$ be an $\N$-gradable $V$-module. We shall use the
same notations as in the preceding section.
Let $T(W)$ be the subspace of elements $w$ of $W$ such that
$u_{n}w=0$ when $\wt u-n-1<0$. The subspace $T(W)$ is called 
the top level of $W$. For $u\in V$, let
$$o(u)=\res_{x}x^{-1}Y_{W}(x^{L(0)}u, x).$$
Then for $w\in T(W)$, $o(u)w\in T(W)$. In fact, it was proved in \cite{Z}
that $u+O(V) \mapsto o(u)$ gives a $A(V)$-module structure on $T(W)$.
Clearly, taking the top level of an $\N$-gradable $V$-module gives 
a functor $T$ from the category of $A(V)$-modules to the category of 
$\N$-gradable $V$-modules.

Consider the affinization $V[t, t^{-1}] = V \otimes \C[t, t^{-1}]$
of $V$ and the tensor algebra $T(V[t, t^{-1}])$ generated by
$V[t, t^{-1}]$. (Note that although we use the same notation,
$T(V[t, t^{-1}])$ has nothing to do with the top level of an $\N$-gradable
weak $V$-module.)  For simplicity, we shall denote $u \otimes t^m$ for
$u \in V$ and $m \in \Z$ by $u(m)$ and we shall omit the tensor product
sign $\otimes$ when we write an element of $T(V[t, t^{-1}])$. Thus
$T(V[t, t^{-1}])$ is spanned by elements of the form $u_1(m_1)\cdots u_k(m_k)$
for $u_i \in V$ and $m_i \in \Z$, $i = 1, \dots, k$.

Let $M$ be an $A(V)$-module and $\rho: A(V) \rightarrow \edo M$
the corresponding representation of $A(V)$.
Consider $T(V[t, t^{-1}])\otimes M$. Again for simplicity we shall
omit the tensor product sign. So $T(V[t, t^{-1}])\otimes M$ is
spanned by elements of the form $u_1(m_1)\cdots u_k(m_k)w$ for
$u_i \in V$, $m_i \in \Z$, $i = 1, \dots, k$ and $w \in M$ and
for any $u \in V$, $m \in Z$, $u(m)$ acts from the left on
$T(V[t, t^{-1}])\otimes M$. For homogeneous $u_i \in V$,
$m_i \in \Z$, $ i = 1, \dots, k$ and $w \in M$, we define the
degree of $u_1(m_1) \cdots u_k(m_k)w$ to be
$(\wt u_1 - m_1 - 1) + \cdots (\wt u_k - m_k - 1)$. In particular,
the degrees of elements of $M$ are $0$.

For $u \in V$, let
\[
Y_t(u, x): T(V[t, t^{-1}]) \otimes M \longrightarrow
T(V[t, t^{-1}]) \otimes M[[x, x^{-1}]]
\]
be defined by
\[
Y_t(u, x) = \sum_{m \in \Z}u(m)x^{-m-1}.
\]
For a homogeneous element $u \in V$, let \[o_t(u) = u(\wt u - 1).\]
Using linearity, we extend $o_t(u)$ to non-homogeneous $u$.

Let $\mathcal{I}$ be the $\Z$-graded
$T(V[t, t^{-1}])$-submodule of
$T(V[t, t^{-1}]) \otimes W$ generated by elements of the forms
$$u(m)w$$
for homogeneous $u \in V$, $w \in T(V[t, t^{-1}])\otimes M$,
$\wt u - m -1+\deg w< 0$,
$$o_t(u)w - \rho(u + O(V))w$$
for $u \in V$, $w \in M$ and
\begin{eqnarray}\label{eq1-1-d}
\lefteqn{\res_{x_{1}}\res_{x_{2}}x_{1}^{p}x_{2}^{q}Y_{t}(u, x_{1})
Y_{t}(v, x_{2})w}\nn
&&-\res_{x_{0}}\res_{x_{2}}p(x_{0}, x_{2})x_{2}^{q}(x_{0}+x_{2})^{\swt u+\deg w}
Y_{t}(Y(u, x_{0})v, x_{2})w,
\end{eqnarray}
where
$$p(x_{0}, x_{2})=\sum_{i=0}^{\swt v+\deg w-q-1}\binom{p-\wt u-\deg w}{ i}
x_{0}^{p-\swt u-\deg w-i}x_{2}^{i},$$
or equivalently,
\begin{eqnarray}\label{eqn2sec3}
u(p)v(q)w &-& \sum_{i = 0}^{\swt v+\deg w-q -1}
\sum_{j = 0}^{\swt u+\deg w}\left(\begin{array}{c}
p - \wt u-\deg w\\ i\end{array}\right)\left(\begin{array}{c}
\wt u+\deg w\\ j\end{array}\right)\cdot\nn
&&\;\;\;\;\;\;\;\;\;\;\cdot (u_{p-\swt u-\deg w-i+j}v)(q+\wt u+\deg w+i-j)w
\end{eqnarray}
for homogeneous $u, v \in V$, $w \in T(V[t, t^{-1}]) \otimes M$,
$p, q \in \Z$  such that $\wt v+\deg w> q$ and
$\wt u - p-1+ \wt v-q - 1 + \deg w \geq 0$.

Let
$S_1(M) = (T(V[t, t^{-1}]) \otimes M)/ \mathcal{I}$. Then $S_1(M)$
is also a $\Z$-graded $T(V[t, t^{-1}])$-module. In fact, by
definition of $\mathcal{I}$, we see that $S_1(M)$ is spanned
by elements of the form $u(m)w + \mathcal{I}$ for homogeneous
$u \in V$, $m \in \Z$ such that $m < \wt u - 1$ and $w \in M$.
In particular, we see that $S_1(M)$ has an $\N$-grading. Note
that since $\mathcal{I} \cap M = \{0\}$, $M$ can be embedded into
$S_1(M)$ and the homogeneous subspace $(S_1(M))_{(0)}$ of degree $0$ is
the image of $M$ under this embedding. We shall identify $M$ with
$(S_1(M))_{(0)}$. In particular, $M$ is now viewed as a subspace of
$S_1(M)$. The vertex operator map $Y_{t}$ induces a vertex operator map
$Y_{S_{1}(M)}: V\otimes S_{1}(M)\to S_{1}(M)[[x, x^{-1}]]$ by
$Y_{S_{1}(M)}(u, x)w=Y_{t}(u, x)w$.

Let $\mathcal{J}$ be the $\N$-graded $T(V[t, t^{-1}])$-submodule
of $S_1(M)$ generated by
\begin{equation}\label{upper-truncation-3-d}
\res_{x_{0}}\res_{x_{2}}x_{0}^{p-\swt u-\deg w-i}x_{2}^{q+i}
(x_{0}+x_{2})^{\swt u+\deg w}
Y_{S_{1}(M)}(Y(u, x_{0})v, x_{2})w.
\end{equation}
for homogeneous $u, v \in V$, $w \in S_1(M)$, $p, q\in \Z$ and
$i\ge \wt v+\deg w-q$.
By Remark \ref{upper-truncation-rmk-2}, $\mathcal{J}$ is generated by
$$\res_{x_{0}}\res_{x_{2}}x_{0}^{m}x_{2}^{K-m-\deg w-2}
(x_{0}+x_{2})^{\deg w}
Y_{S_{1}(M)}(Y((x_{0}+x_{2})^{L(0)}u, x_{0})x_{2}^{L(0)}v, x_{2})
w$$
for homogeneous $u, v \in V$, $w \in S_1(M)$, $K, m\in \Z$ satisfying
$m\le K-2\deg w-2$,
or equivalently, by
\begin{equation}
\sum_{j=0}^{\swt u+\deg w}\binom{\wt u+\deg w}{j}(u_{j+m}v)(K+\wt u+\wt v- j - m -2)w,
\end{equation}
for homogeneous $u, v \in V$ and the same $w$, $K$ and $m$.

Let $S(M) = S_1(M)/\mathcal{J}$. Then $S(M)$
is also an $\N$-graded $T(V[t, t^{-1}])$-module. We can still
use elements of $T(V[t, t^{-1}])\otimes M$ to represent elements
of $S(M)$. But note that these elements now satisfy more relations.
We equip $S(M)$ with the vertex operator map
\[
Y_{S(M)}: V \otimes S(M) \longrightarrow S(M)[[x, x^{-1}]]
\]
given by
\[
u \otimes w \rightarrow Y_{S(M)}(u, x)w = Y_{S_{1}(M)}(u, x)w=Y_{t}(u, x)w.
\]

\begin{thm}\label{functor-1}
The pair
$(S(M), Y_{S(M)})$ is an $\N$-gradable weak $V$-module.
\end{thm}
\pf
As in $S_1(M)$, for $u \in V$ and $w \in S(M)$, we
have $u(m)w = 0$ when $m > \wt u + \wt w - 1$. Clearly,
\[
Y({\bf 1}, x) = 1_{S(M)},
\]
where $1_{S(M)}$ is the identity operator on $S(M)$.

By Theorems \ref{main proposition} and \ref{main theorem},
$S(M)$ is an $\N$-gradable weak $V$-module. \epfv

Let $M_{1}$ and $M_{2}$ be $\N$-gradable weak $V$-modules and
$f: M_{1} \to M_{2}$ a module map. Then $F$ induces a linear map
from $T(V[t, t^{-1}])\otimes M_{1}$ to $T(V[t, t^{-1}])\otimes M_{2}$.
By definition, this induced linear map in turn induces a linear map
$S(f)$ from $S(M_{1})$ to $S(M_{2})$. Since $Y_{S(M_{1})}$ and $Y_{S(M_{2})}$
are induced by $Y_{t}$ on $T(V[t, t^{-1}]) \otimes M_{1}$ and
$T(V[t, t^{-1}]) \otimes M_{2}$, respectively, we have
$$S(f)(Y_{S(M_{1})}(u, x)w_{1})= Y_{S(M_{2})}(u, x)S(f)(w_{1})$$
for $u\in V$ and $w_{1}\in S(M_{1})$. Thus $S(f)$ is a module map.
The following result is now clear:

\begin{cor}\label{functor-2}
Let $V$ be a $\N$-graded vertex algebra.
Then the correspondence sending an $A(V)$-module $M$ to an $\N$-gradable
weak $V$-module $(S(M), Y_{S(M)})$ and an $A(V)$-module map $M_{1}\to M_{2}$
to a $V$-module map $S(f): S(M_{1})\to S(M_{2})$ is a functor from the category of
$A(V)$-modules to the category of $\N$-gradable
weak $V$-modules.
\end{cor}

\section{A universal property and $T$ as a left inverse of $S$}

In this section, we prove that $S$ satisfies a natural universal
property and thus is the same as the functor constructed in
\cite{DLM}. In particular, we achieve our goal of constructing
$\N$-gradable weak $V$-modules from $A(V)$-modules without dividing
relations corresponding to the commutator formula for weak modules.
We also prove that the functor $T$ of sending an $\N$-gradable weak $V$-module
$W$ to its top level $T(W)$ is a left inverse of $S$.

As in the preceding section, we let $V$ be a $\N$-graded vertex algebra.

We shall need the following lemma:

\begin{lemma}\label{lemma1sec3}
In the setting of Theorem \ref{main proposition}, if we assume only
(\ref{eq1-1}) (or equivalently, (\ref{eq1})), then
for $u, v\in V$, homogeneous $w\in W$ and
$p', q'\in \Z$, $q'< \deg w$, we have
\begin{eqnarray}\label{eq1.1sec3}
\lefteqn{\res_{x_1}\res_{x_2}x_{1}^{p'}x_{2}^{q'}Y_{W}(x_{1}^{L(0)}u, x_{1})
Y_{W}(x_{2}^{L(0)}v, x_{2})w }\nn
&&= \res_{x_0}\res_{x_2}q(x_{0}, x_{2})x_2^{q'}(x_0 + x_2)^{\deg w}
Y_W(Y((x_0 + x_2)^{L(0)}u, x_0)x_{2}^{L(0)}v, x_2)w,
\end{eqnarray}
where $q(x_{0}, x_{2})\in C[x_{0}, x_{0}^{-1}, x_{2}]$
such that when $w$ is homogeneous,
$$q(x_{0}, x_{2})=\sum_{i=0}^{\deg w - q'-1}
\binom{p'-\deg w}{i}
x_{0}^{p'-\deg w-i}x_{2}^{i}.$$
\end{lemma}
\pf In the case that $u, v$ are also homogeneous,
the lemma follows from Remark \ref{eq1-rmk} by
taking $p=p'+\wt u$, $q=q'+\wt v$, $k=\wt v+\deg w$ and
$l=\wt u+\deg w$. \epfv

For an $A(V)$-module $M$, let $\mathcal{I}$, $S_{1}(M)$ and $\mathcal{J}$
be the same as in the preceding section. The following theorem is the
main technical result in this section
and its proof is the hardest in the present paper:

\begin{thm}\label{lemma2sec3}
Let $M$ be a an $A(V)$-module. Then in $S_1(M)$,
\[
\mathcal{J} \cap M = 0.
\]
\end{thm}
\pf
Note that elements of the form (\ref{eqn2sec3}) are in
$\mathcal{I}$ and that
$M$ is identified with $(S_1(M))_{(0)}$. Using the definitions of $\mathcal{I}$
and $\mathcal{J}$ and noticing that the degree of elements of $M$ is $0$,
we see that elements of
$\mathcal{J}$ are linear combinations of
elements of the form
\begin{eqnarray}\label{eqn3-1sec3}
\lefteqn{a(\wt a-K+N-1)\res_{x_{0}}\res_{x_{2}}x_{0}^{m}x_{2}^{K-m-N-2}\cdot}\nn
&&\quad\quad\quad\quad\quad\cdot (x_{0}+x_{2})^{N}
Y_{S_{1}(M)}(Y((x_{0}+x_{2})^{L(0)}u, x_{0})x_{2}^{L(0)}v, x_{2})b(\wt b-N-1)w+\mathcal{I}\nn
&&=a(\wt a-K+N-1)\sum_{j=0}^{\swt u+N}\binom{\wt u+N}{ j}\cdot\nn
&&\quad\quad\quad\quad\quad\cdot
(u_{j+m}v)(K+\wt u+\wt v- j - m -2)b(\wt b-N-1)w
+\mathcal{I}\nn
\end{eqnarray}
for homogeneous $a, b, u, v \in V$, $w \in M$, $K, N, m\in \Z$
satisfying $m\le K-2N-2$. To prove this theorem, we need only prove that
(\ref{eqn3-1sec3}) is $0$ in $S_1(M)$.

When $N<0$, $\deg b(\wt b-N-1)w=N<0$. So $b(\wt b-N-1)w\in \mathcal{I}$
and hence (\ref{eqn3-1sec3}) is $0$ in $S_1(M)$.

When $N\ge 0$ but $K>N$,
$$\deg (u_{j+m}v)(K+\wt u+\wt v- j - m -2)b(\wt b-N-1)w=N-K<0.$$
So $(u_{j+m}v)(K+\wt u+\wt v- j - m -2)b(\wt b-N-1)w\in \mathcal{I}$
and hence (\ref{eqn3-1sec3}) is $0$ in $S_1(M)$.

We now prove that (\ref{eqn3-1sec3}) is $0$ in $S_1(M)$ in the case
$N\ge 0$ and $N\ge K$.
We rewrite (\ref{eqn3-1sec3}) as
\begin{eqnarray}\label{eqn3-2sec3}
\lefteqn{\res_{y_{1}}\res_{y_{2}}\res_{x_{0}}\res_{x_{2}}
y_{1}^{-(K-N)-1}y_{2}^{-N-1}x_{0}^{m}x_{2}^{K-m-N-2}
(x_{0}+x_{2})^{N}\cdot}\nn
&&\quad\quad\quad \quad\cdot
Y_{S_{1}(M)}(y_{1}^{L(0)}a, y_{1})
Y_{S_{1}(M)}(Y((x_{0}+x_{2})^{L(0)}u, x_{0})x_{2}^{L(0)}v, x_{2})
Y_{S_{1}(M)}(y_{2}^{L(0)}b, y_{2})w\nn
\end{eqnarray}
where $a, b, u, v \in V$, $w \in M$, $K, N, m\in \Z$
satisfying $m\le K-2N-2$.

Note that in $S_{1}(M)$, $Y_{S_{1}(M)}$ satisfies (\ref{eq1.1sec3}) with $Y_{W}$
replaced by $Y_{S_{1}(M)}$. We shall use (\ref{eq1.1sec3}) and other formulas
that still hold to show that (\ref{eqn3-2sec3}) is in fact
a linear combination of an element of $O(V)$ acting on $w \in M$
and hence equals $0$ in $S_1(M)$.

Let
$$\tilde{w}=\res_{y_{2}}y_{2}^{-N-1}Y_{S_{1}(W)}(y_{2}^{L(0)}b, y_{2})w.$$
Using the properties of vertex operators, changing variables and using
(\ref{eq1.1sec3}), we see that  (\ref{eqn3-2sec3}) is equal to
\begin{eqnarray}\label{eqn3-3sec3}
\lefteqn{\res_{y_{1}}\res_{x_{0}}\res_{x_{2}}
y_{1}^{-(K-N)-1}x_{0}^{m}x_{2}^{K-m-N-2}
(x_{0}+x_{2})^{N}\cdot}\nn
&&\quad\quad\quad \quad\cdot
Y_{S_{1}(M)}(y_{1}^{L(0)}a, y_{1})
Y_{S_{1}(M)}(Y((x_{0}+x_{2})^{L(0)}u, x_{0})x_{2}^{L(0)}v, x_{2})\tilde{w}\nn
&&=\res_{y_{1}}\res_{x_{0}}\res_{x_{2}}
y_{1}^{-(K-N)-1}x_{0}^{m}x_{2}^{K-m-2}
(x_{0}x_{2}^{-1}+1)^{N}\cdot\nn
&&\quad\quad\quad \quad\cdot
Y_{S_{1}(M)}(y_{1}^{L(0)}a, y_{1})
Y_{S_{1}(M)}(x_{2}^{L(0)}Y((x_{0}x_{2}^{-1}+1)^{L(0)}u, x_{0}x_{2}^{-1})v, x_{2})\tilde{w}\nn
&&=\res_{x_{0}}x_{0}^{m}(x_{0}+1)^{N}\res_{y_{1}}\res_{x_{2}}
y_{1}^{-(K-N)-1}x_{2}^{K-1}
\cdot\nn
&&\quad\quad\quad \quad\cdot
Y_{S_{1}(M)}(y_{1}^{L(0)}a, y_{1})
Y_{S_{1}(M)}(x_{2}^{L(0)}Y((x_{0}+1)^{L(0)}u, x_{0})v, x_{2})\tilde{w}\nn
&&=\res_{x_{0}}x_{0}^{m}(x_{0}+1)^{N}\res_{y_{0}}\res_{x_{2}}
q_{1}(y_{0}, x_{2})x_{2}^{K-1}(y_{0}+x_{2})^{N}
\cdot\nn
&&\quad\quad\quad \quad\cdot
Y_{S_{1}(M)}(Y((y_{0}+x_{2})^{L(0)}a, y_{0})
x_{2}^{L(0)}Y((x_{0}+1)^{L(0)}u, x_{0})v, x_{2})\tilde{w},
\end{eqnarray}
where
\begin{eqnarray*}
q_{1}(y_{0}, x_{2})&=&\sum_{i=0}^{N - K}
\binom{-K-1}{ i}
y_{0}^{-K-1-i}x_{2}^{i}\nn
&=&x_{2}^{-K-1}\sum_{i=0}^{N - K}
\binom{-K-1}{ i}
(y_{0}x_{2}^{-1})^{-K-1-i}\nn
&=&x_{2}^{-K-1}q_{1}(y_{0} x_{2}^{-1}, 1).
\end{eqnarray*}

The right-hand side of (\ref{eqn3-3sec3}) can now be written as
\begin{eqnarray}\label{eqn3-4sec3}
\lefteqn{\res_{x_{0}}x_{0}^{m}(x_{0}+1)^{N}\res_{y_{0}}\res_{x_{2}}
q_{1}(y_{0} x_{2}^{-1}, 1)x_{2}^{N-2}(y_{0}x_{2}^{-1}+1)^{N}
\cdot}\nn
&&\quad\quad\quad \quad\cdot
Y_{S_{1}(M)}(x_{2}^{L(0)}Y((y_{0}x_{2}^{-1}+1)^{L(0)}a, y_{0}x_{2}^{-1})
Y((x_{0}+1)^{L(0)}u, x_{0})v, x_{2})\tilde{w}\nn
&&=\res_{x_{0}}x_{0}^{m}(x_{0}+1)^{N}\res_{y_{0}}\res_{x_{2}}
q_{1}(y_{0}, 1)x_{2}^{N-1}(y_{0}+1)^{N}
\cdot\nn
&&\quad\quad\quad \quad\cdot
Y_{S_{1}(M)}(x_{2}^{L(0)}Y((y_{0}+1)^{L(0)}a, y_{0})
Y((x_{0}+1)^{L(0)}u, x_{0})v, x_{2})\tilde{w}.
\end{eqnarray}

Using the commutator formula for $V$, the right-hand side of (\ref{eqn3-4sec3})
is equal to
\begin{eqnarray}\label{eqn3-5sec3}
\lefteqn{\res_{y_{0}}q_{1}(y_{0}, 1)(y_{0}+1)^{N}\res_{x_{2}}\res_{x_{0}}
x_{2}^{N-1}x_{0}^{m}(x_{0}+1)^{N}
\cdot}\nn
&&\quad\quad\quad \quad\cdot
Y_{S_{1}(M)}(x_{2}^{L(0)}
Y((x_{0}+1)^{L(0)}u, x_{0})Y((y_{0}+1)^{L(0)}a, y_{0})v, x_{2})\tilde{w}\nn
&&\quad +\res_{y_{0}}q_{1}(y_{0}, 1)(y_{0}+1)^{N}\res_{x_{2}}\res_{x_{0}}\res_{x}
x_{2}^{N-1}x_{0}^{m}(x_{0}+1)^{N}\cdot\nn
&&\quad\quad\quad \quad\cdot
x_{0}^{-1}\left(\frac{y_{0}-x}{x_{0}}\right)
Y_{S_{1}(M)}(x_{2}^{L(0)}Y(Y((y_{0}+1)^{L(0)}a, x)(x_{0}+1)^{L(0)}u, x_{0})v, x_{2})\tilde{w}.\nn
\end{eqnarray}

The first term of (\ref{eqn3-5sec3}) is a linear combination of
the elements of the form
\begin{equation}\label{eqn3-6sec3}
\res_{x_{2}}\res_{x_{0}}x_{2}^{N-1}x_{0}^{n}(x_{0}+1)^{N}
Y_{S_{1}(M)}(x_{2}^{L(0)}
Y((x_{0}+1)^{L(0)}\tilde{u}, x_{0})\tilde{v}, x_{2})\tilde{w}
\end{equation}
for $n=m\le -N-2$, $\tilde{u}=u, \tilde{v}\in V$.
The second term of (\ref{eqn3-5sec3}) is equal to
\begin{eqnarray}\label{eqn3-7sec3}
\lefteqn{\res_{x_{2}}\res_{x_{0}}\res_{x}\res_{y_{0}}q_{1}(x_{0}+x, 1)(x_{0}+x+1)^{N}
x_{2}^{N-1}x_{0}^{m}(x_{0}+1)^{N}x_{0}^{-1}\left(\frac{y_{0}-x}{x_{0}}\right)\cdot}\nn
&&\quad\quad\cdot
Y_{S_{1}(M)}(x_{2}^{L(0)}Y((x_{0}+1)^{L(0)}
Y((1+x(x_{0}+1)^{-1})^{L(0)}a, x(x_{0}+1)^{-1})u, x_{0})v, x_{2})\tilde{w}\nn
&&=\res_{x_{2}}\res_{x_{0}}\res_{x}q_{1}(x_{0}+x, 1)(x_{0}+x+1)^{N}
x_{2}^{N-1}x_{0}^{m}(x_{0}+1)^{N}\nn
&&\quad\quad\cdot
Y_{S_{1}(M)}(x_{2}^{L(0)}Y((x_{0}+1)^{L(0)}
Y((1+x(x_{0}+1)^{-1})^{L(0)}a, x(x_{0}+1)^{-1})u, x_{0})v, x_{2})\tilde{w}\nn
&&=\res_{x_{2}}\res_{x_{0}}\res_{y}q_{1}(x_{0}+y(x_{0}+1), 1)(x_{0}+1)^{N}(1+y)^{N}
x_{2}^{N-1}x_{0}^{m}(x_{0}+1)^{N+1}\nn
&&\quad\quad\cdot
Y_{S_{1}(M)}(x_{2}^{L(0)}Y((x_{0}+1)^{L(0)}
Y((1+y)^{L(0)}a, y)u, x_{0})v, x_{2})\tilde{w}.
\end{eqnarray}
By the definition of $q_{1}(y_{0}, x_{2})$ above, we see
that $q_{1}(x_{0}+y(x_{0}+1), 1)$ is a Laurent series in $x_{0}$ and $y$
such that the largest power of $x_{0}$ in $q_{1}(x_{0}+y(x_{0}+1), 1)$
is $-K-1$. Then the powers
of $x_{0}$ in
$$q_{1}(x_{0}+y(x_{0}+1), 1)
x_{0}^{m}(x_{0}+1)^{N+1}$$
are less than or equal to $-N-2$. Thus the right-hand side of
(\ref{eqn3-7sec3}) is a linear
combination of the elements of the form (\ref{eqn3-6sec3})
for $n\le -N-2$ and $\tilde{u}, \tilde{v}=v\in V$.

It remains to prove that (\ref{eqn3-6sec3}) for $n\le -N-2$ and $\tilde{u},
\tilde{v}\in V$
is $0$ in $S_{1}(M)$.
We prove that elements of the form (\ref{eqn3-6sec3}) for $n\le -N-2$ and $\tilde{u},
\tilde{v}\in V$
are linear
combinations of elements obtained from the action of
$O(V)$ on $M$. Then by the definition of $S_{1}(M)$,
(\ref{eqn3-6sec3}) for $n\le -N-2$ and $\tilde{u}, \tilde{v}\in V$ is $0$.

Writing the element $\tilde{w}$ explicitly and using (\ref{eq1.1sec3}),
we see that (\ref{eqn3-6sec3})
is equal to
\begin{eqnarray}\label{eqn3-8sec3}
\lefteqn{\res_{x_{0}}x_{0}^{n}(x_{0}+1)^{N}
\res_{x_{2}}\res_{y_{2}}x_{2}^{N-1}y_{2}^{-N-1}
Y_{S_{1}(M)}(x_{2}^{L(0)}
Y((x_{0}+1)^{L(0)}\tilde{u}, x_{0})\tilde{v}, x_{2})
Y_{S_{1}(M)}(y_{2}^{L(0)}b, y_{2})w}\nn
&&=\res_{x_{0}}x_{0}^{n}(x_{0}+1)^{N}\res_{y_{0}}\res_{y_{2}}q_{2}(y_{0}, y_{2})
y_{2}^{-N-1}\cdot\nn
&&\quad\quad\quad\quad\quad\quad\cdot
Y_{S_{1}(M)}(Y((y_{0}+y_{2})^{L(0)}
Y((x_{0}+1)^{L(0)}\tilde{u}, x_{0})\tilde{v}, y_{0})
y_{2}^{L(0)}b, y_{2})w\nn
&&=\res_{x_{0}}x_{0}^{n}(x_{0}+1)^{N}\res_{y_{0}}\res_{y_{2}}q_{2}(y_{0}, y_{2})
y_{2}^{-N-1}\cdot\nn
&&\quad\quad\quad\quad\quad\quad\cdot
Y_{S_{1}(M)}(y_{2}^{L(0)}Y((y_{0}y_{2}^{-1}+1)^{L(0)}
Y((x_{0}+1)^{L(0)}\tilde{u}, x_{0})\tilde{v}, y_{0}y_{2}^{-1})
b, y_{2})w\nn
&&=\res_{x_{0}}x_{0}^{n}(x_{0}+1)^{N}\res_{y}\res_{y_{2}}q_{2}(yy_{2}, y_{2})
y_{2}^{-N}\cdot\nn
&&\quad\quad\quad\quad\quad\quad\cdot
Y_{S_{1}(M)}(y_{2}^{L(0)}Y((y+1)^{L(0)}
Y((x_{0}+1)^{L(0)}\tilde{u}, x_{0})\tilde{v}, y)
b, y_{2})w\nn
&&=\res_{x_{0}}x_{0}^{n}(x_{0}+1)^{N}\res_{y}\res_{y_{2}}q_{2}(yy_{2}, y_{2})
y_{2}^{-N}\cdot\nn
&&\quad\quad\quad\quad\quad\quad\cdot
Y_{S_{1}(M)}(y_{2}^{L(0)}Y(
Y((x_{0}+1)^{L(0)}(y+1)^{L(0)}\tilde{u}, x_{0}(y+1))(y+1)^{L(0)}\tilde{v}, y)
b, y_{2})w\nn
&&=\res_{y_{2}}\res_{y}\res_{x}x^{n}(x+y+1)^{N}q_{2}(yy_{2}, y_{2})
y_{2}^{-N}(y+1)^{-N-1-n}\cdot\nn
&&\quad\quad\quad\quad\quad\quad\cdot
Y_{S_{1}(M)}(y_{2}^{L(0)}Y(
Y((x+y+1)^{L(0)}\tilde{u}, x)(y+1)^{L(0)}\tilde{v}, y)
b, y_{2})w\nn
&&=\res_{y_{2}}\res_{y}\res_{x_{1}}(x_{1}-y)^{n}(x_{1}+1)^{N}q_{2}(yy_{2}, y_{2})
y_{2}^{-N}(y+1)^{-N-1-n}\cdot\nn
&&\quad\quad\quad\quad\quad\quad\cdot
Y_{S_{1}(M)}(y_{2}^{L(0)}Y((x_{1}+1)^{L(0)}\tilde{u}, x_{1})Y((y+1)^{L(0)}\tilde{v}, y)
b, y_{2})w\nn
&&\quad -\res_{y_{2}}\res_{y}\res_{x_{1}}(-y+x_{1})^{n}(x_{1}+1)^{N}q_{2}(yy_{2}, y_{2})
y_{2}^{-N}(y+1)^{-N-1-n}\cdot\nn
&&\quad\quad\quad\quad\quad\quad\cdot
Y_{S_{1}(M)}(y_{2}^{L(0)}Y((y+1)^{L(0)}\tilde{v}, y)
Y((x_{1}+1)^{L(0)}\tilde{u}, x_{1})
b, y_{2})w\nn
\end{eqnarray}
Where
$$q_{2}(y_{0}, y_{2})=\sum_{i=0}^{N}\binom{N-1}{ i}y_{0}^{N-1-i}y_{2}^{i}.$$

From the expression of $q_{2}(y_{0}, y_{2})$, we have
$$q_{2}(yy_{2}, y_{2})=y_{2}^{N-1}q_{2}(y, 1).$$
Thus the first term in the right-hand side of (\ref{eqn3-8sec3})
is a linear combination of
terms of the form
$$\res_{y_{2}}y_{2}^{-1}
Y_{S_{1}(M)}(y_{2}^{L(0)}\res_{x_{1}}x_{1}^{\tilde{n}}
Y((x_{1}+1)^{L(0)}\tilde{u}, x_{1})\tilde{\tilde{v}}, y_{2})w$$
for $\tilde{n}\le -2$ and $\tilde{\tilde{v}}\in V$.
By definition,
$$\tilde{\tilde{u}}=\res_{x_{1}}x_{1}^{\tilde{n}}
Y((x_{1}+1)^{L(0)}\tilde{u}, x_{1})\tilde{\tilde{v}}\in O(V)$$
and, by definition,
$$\res_{y_{2}}y_{2}^{-1}
Y_{S_{1}(M)}(y_{2}^{L(0)}\tilde{\tilde{u}}, y_{2})=o(\tilde{u}).$$
Since $o(\tilde{\tilde{u}})w=0$ in $S_{1}(M)$,
$\res_{y_{2}}y_{2}^{-1}
Y_{S_{1}(M)}(y_{2}^{L(0)}\tilde{\tilde{u}}, y_{2})w$ is also $0$. This proves that
the first term in the right-hand side of (\ref{eqn3-8sec3}) is $0$.

The largest power of $y$ in
$$(-y+x_{1})^{n}q_{2}(y, 1)(y+1)^{-N-1-n}$$
is $-2$. Thus the second term in the right-hand side of (\ref{eqn3-8sec3})
is a linear combination of
terms of the form
$$\res_{y_{2}}y_{2}^{-1}
Y_{S_{1}(M)}(y_{2}^{L(0)}\res_{y}y^{\tilde{n}}Y((y+1)^{L(0)}\tilde{v}, y)
\tilde{\tilde{u}}, y_{2})w$$
for $\tilde{n}\le -2$ and $\tilde{\tilde{u}}\in V$. Then the same argument
as above shows that the second term in the right-hand side of (\ref{eqn3-8sec3})
is $0$. \epfv

\begin{rema}
{\rm Note that though the proof of the theorem above does not
need the commutator formula for $Y_{S_{1}(M)}$, the commutator formula
for $V$ is indeed used. It will be interesting to see whether
there is a proof without using the commutator formula for $V$.}
\end{rema}

\begin{cor}\label{embedding}
Let $M$ be an $A(V)$-module. Then the embedding from $M$ to $S_{1}(M)$
induces an isomorphism $e_{M}$ of $A(V)$-modules from $M$ to $T(S(M))$.
\end{cor}
\pf
Given an $N$-gradable weak $V$-module $M$, we have an embedding from
$M$ to $S_{1}(M)$. By Theorem \ref{lemma2sec3}, the induced map from
$M$ to $S(M)=S_{1}/\mathcal{J}$ is still an embedding. By definition,
$T(S(M))$ is exactly the image of $M$ in $S(M)$. Clearly,
when we view the induced map from
$M$ to $S(M)=S_{1}/\mathcal{J}$ as a map from $M$ to $T(S(M))$,
it is an $A(V)$-module map and thus is a
an isomorphism of $A(V)$-modules.
\epfv

\begin{thm}
The functor $S$ has the following universal property:
Let $M$ be an $A(V)$-module. For any $\N$-gradable weak $V$-module
$W$ and any
$A(V)$-module map $f: M\to T(W)$, there exists a unique
$V$-module map $\tilde{f}: S(M)\to W$ such that
$\tilde{f}|_{T(S(M)}=f\circ e_{M}^{-1}$.
\end{thm}
\pf
Elements of $S(M)$ are linear combinations of elements of the form $u(m)w$
for homogeneous $u\in V$, $w\in M$ and $m\in \Z$ satisfying
$\wt -m-1+\deg w\ge 0$. We define $\tilde{f}(u(m)w)=u_{m}f(w)$.
By definition, relations among elements of the form $u(m)w$
for homogeneous $u\in V$, $w\in M$ and $m\in \Z$ satisfying
$\wt -m-1+\deg w\ge 0$ are also satisfied by elements of $W$ of the form
$u_{m}f(w)$. Thus the map $\tilde{f}$ is well defined. Using
(\ref{eq1-0}) (or equivalently,
(\ref{eq1}) or (\ref{eq1-1})), we see that $\tilde{f}$ is indeed a
$V$-module map. By definition, $\tilde{f}|_{M}=f\circ e_{M}^{-1}$. If $g: S(M)\to W$
is another such $V$-module map,
then $g(u(m)w)=u_{m}g(w)=u_{m}f(w)=\tilde{f}(u(m)w)$, proving the uniqueness of
$\tilde{f}$. \epfv

We also have the following:

\begin{cor}
The functor $T$ is a left inverse of $S$, that is, $T\circ S$ is naturally 
isomorphic to is the identity functor on the category of $A(V)$-modules.
\end{cor}
\pf
For an $A(V)$-module $M$, by Corollary \ref{embedding}, 
we have an isomorphism $e_{M}^{-1}$ from $(T\circ S)(M)$ to $M$. 
For an $A(V)$-module map
$f: M_{1}\to M_{2}$, by definition, the 
$A(V)$-module map $(T\circ S)(f): (T\circ S)(M_{1})\to (T\circ S)(M_{2})$ 
is $e_{M_{2}}\circ f\circ e_{M_{1}}^{-1}$. Thus we obtain a natural 
transformation $e^{-1}$ 
from $T\circ S$ to the identity functor on the category of $A(V)$-modules.
This natural transformation $e^{-1}$ is clearly a natural isomorphism 
since it has an inverse $e$ which is given by the isomorphism $e_{M}$
from an $A(V)$-module $M$ to $T(S(M))$.
\epfv

In \cite{Z}, in his proof that the set of equivalence classes 
of irreducible $\N$-gradable weak $V$-modules is in bijection with
the set of equivalence classes of irreducible modules for $A(V)$, Zhu
in fact obtained implicitly a functor from the category of $A(V)$-modules
to the category of $\N$-gradable weak $V$-modules. 
In \cite{DLM}, this functor, denoted by $\bar{M}_{0}$,
was constructed explicitly and was proved to satisfy the same universal property 
above. The following result achieves our goal
of constructing this functor without dividing
relations corresponding to the commutator formula for weak modules:

\begin{cor}
The functor $S$ is equal to
the functor $\bar{M}_{0}$ constructed in \cite{DLM}.
\end{cor}
\pf 
This result follows immediately from the universal property.
\epfv

\def\refname{\hfil{REFERENCES}}

\end{document}